\documentclass[12pt]{amsart}
\usepackage[all, cmtip]{xy}
\usepackage{amscd, amssymb, amsmath, latexsym}

\def\qed{{\hfill \bf q.e.d.}}

\newcommand{\bq}{\begin{equation}}
\newcommand{\eq}{\end{equation}}
\newcommand{\ba}{\begin{array}}
\newcommand{\ea}{\end{array}}
\def \bc {\begin{center}}
\def \ec {\end{center}}

\newtheorem{thm}{Theorem}[section]

\newtheorem{lemma}[thm]{Lemma}

\newtheorem{conj}[thm]{Conjecture}

\begin{document}

\title[Half-isomorphisms]{Half-isomorphisms of finite automorphic
Moufang loops.} 

\author{A. Grishkov}

\address{A.Grishkov, Instituto de Matem\'atica e Estat\'istica\\Universidade de S\~ao Paulo\\
 Rua de Mat\~ao 1010\\S\~ao Paulo, 05508-090\\
Brazil}
\email{grishkov@ime.usp.br}

\author{M.L.Merlini Giuliani}

\address{M.L.Merlini Giuliani, Universidade Federal do ABC, Av. dos Estados, 5001 – Bairro Bangu, Santo Andr\'e – SP 09210-580 - Brazil\\}

\email{maria.giuliani@ufabc.edu.br}

\author{M. Rasskazova}

\address{M. Rasskazova, Omsk State Institute of Service\\
Pevtsova street 13\\ Omsk, 644099\\ Russia}
\email{marinarasskazova@yandex.ru}

\author{L. Sabinina}

\address{L. Sabinina, Universidad Aut\'onoma del Estado de Morelos\\
Av. Universidad 1001\\ Cuernavaca, 62209\\ M\'exico}
\email{liudmila@uaem.mx}

\begin{abstract}

We show that each half-automorphism of a finite automorphic Moufang loop is trivial. In general this is not true for finite left automorphic Moufang 
loops and for finite automorphic loops.
\end{abstract}

\maketitle

\section{Introduction.}

A half-isomorphism of a loop is bijection $\tau$  such that $\tau(xy)$ equals either $\tau(x)\tau(y)$ or  $\tau(y)\tau(x)$ for
all elements $x,y$.  For groups all half-isomorphisms are either an isomorphism or an anti-isomorphism as was shown by Scott \cite{Scott}.
Following the terminology of \cite{Scott} we will call 
 such half-isomorphisms trivial. In general for loops the situation is different: there exist Moufang loops with half-isomorphisms which are neither an isomorphism nor an anti-isomorphism.
In many cases, in particular, for finite Moufang loops of odd order, all half-isomorphisms are trivial  (\cite{GG},\cite{L}).   Obviously the half-isomorphisms of a loop always form a group, which is a new and, at the time of the writing this note, rather mysterious invariant of a loop.
Let us note that any anti-isomorphism of a Moufang loop $M$ to a Mou\-fang loop $M_1$ is the composition of an
isomorphism of $M$ to $M_1$ and inverse map on  $M_1$.

\section{Preliminaries and definitions.}
In this paper we consider a
 {\it loop} as a universal algebra $\langle Q;\cdot, \backslash, /,1\rangle$ of type 
$(3,0,1)$ such that the identities

$$(x\cdot y)/y=x = (x/y)\cdot y,$$ $$ x\cdot (x\backslash y) =y=x\backslash (x\cdot y), $$
$$x\cdot 1=x=1\cdot x$$
hold for all $x,y\in L$.
In the following we often write $xy$ instead of $x\cdot y$.\\

A {\it Moufang loop}
is a loop in which one of  the following equivalent identities holds:
$$((xy)x)z=x(y(xz)),$$
$$((xy)z)y =x(y(zy)),$$
$$(xy)(zx)=(x(yz))x.$$\\
The bijection  $L_a: Q\to Q,\,\, L_ax=ax$ is called  a {\it left multiplication}, analogously the  bijection $R_a:Q\to Q, \,\,R_ay=ya$ is called a {\it right multiplication} for all  $\,a\in Q$. 

With an arbitrary loop $Q$ one may associate several groups of transformations of the set $Q$, such as:
\begin{itemize}
\item the multiplication group $\mathrm{Mlt}(Q)$ generated by the left and right multiplications by elements of $Q$;
\item the left multiplication group $\mathrm{LMlt}(Q)\subseteq \mathrm{Mlt}(Q)$ generated by the left multiplications only;
\item the inner mapping group $\mathrm{Inn}(Q)\subseteq \mathrm{Mlt(Q)}$ defined as the stabilizer of the neutral element $1\in Q$;
\item the left inner mapping group $\mathrm{LInn}(Q)$ which is the stabilizer of $1$ in  $\mathrm{LMlt}(Q)$.
\end{itemize}

$\mathrm{Inn}(Q),$ the  {\it inner mapping group} of $Q$,  is 
generated by three families of elements of  $ \mathrm{Mlt(Q)}$:
$$\ell_{x,y}= L^{-1}_{xy}\circ L_x\circ L_y,$$
$$r_{x,y}= R^{-1}_{xy}\circ R_y\circ R_x,$$
$$T_x= L^{-1}_{x}\circ R_x$$
for all $x,y\in Q.$
The loop $Q$ is called {\it automorphic}  if $\mathrm{Inn}(Q)$ acts on $Q$ by automorphisms, and { \it  left automorphic}  if $\mathrm{LInn}(Q)$ does.

A loop $Q$ is  automorphic if the  mappings $\ell_{x,y},\,\, r_{x,y}\,$ and $\,T_x$ for all $x,y\in Q$ are automorphisms of $Q$. 
A loop $Q$ is  left automorphic
if the mappings $\ell_{x,y}$ are automorphisms of $Q$. It is known  for a Moufang loop $Q$
that if $\ell_{x,y}$ is an automorphism for some $x,y \in Q$ then $\,r_{x,y}$ is also an automorphism of $Q$ and vice versa (see \cite {Bruck}).

The property of a loop to be left automorphic is of significance in differential geometry. In particular, a reductive homogeneous 
space is completely determined by a loop
with the property of left power-alternativity (or left monoalternativity) and the left automorphic property (\cite{Kikkawa}, \cite{Sabinin}).

Each alternative (in particular, each diassociative) automorphic loop
is Moufang (\cite {KKP}).
Also, each power-alternative left automorphic  smooth loop is
Moufang (\cite {CS}).
Besides every smooth automorphic Moufang loop is a group. A well known
class of non-trivial finite examples of left automorphic Moufang loops is
provided by Code loops  (\cite {Griess}).

As usual one defines  a {\it normal subloop} as the kernel of a loop homomorphism.
A subloop of a given loop $Q$ is normal if and only if it is invariant under
the $\mathrm{Inn}(Q)$. 

Moufang loops are diassociative, that is the subloop generated by two elements is a group.
Hence in a  Moufang loop one has $$x\backslash y=x^{-1}y,\,\,y/x=yx^{-1}.$$

In a Moufang loop $Q$ the {\it commutator subloop} $[Q,Q]$ is  the  subloop generated
by all elements of the form 
$x^{-1}y^{-1}xy= [x,y]$, where $x,y \in Q$ and
the  {\it associator subloop}  $(Q,Q,Q)$ is the subloop generated by all elements of the form $(x(yz)) ^{-1}((xy)z)=(x,y,z)$ for all $x,y,z \in Q.$

The subloop  $Z(Q)=\{x\in Q:[x,a]=1,\forall a\in Q\}$ is called the  {\it commutant} of $Q.$ 

The subloops 
$$N_\lambda (Q)=\{x\in Q:(x,a,b)=1,\forall a,b\in Q\},$$
$$N_\mu (Q) =\{x\in Q:(a,x,b)=1,\forall a,b\in Q\},$$
$$N_\rho (Q)=\{x\in Q:(a,b,x)=1,\forall a,b\in Q\}$$
 are  called {\it left, middle} and  {\it right   nucleus} of $Q$. 
The subloop $N(Q) = N_\lambda (Q)\cap N_\mu (Q)\cap N_\rho (Q)$ 
is called the {\it  nucleus} of $Q$.
In a Moufang loop $L$ all nuclei coincide.
The subloop 
$C(Q)=Z(Q)\cap N(Q)$ is called  the {\it center} of the loop $Q$.

All characteristic subloops of any automorphic loop are normal. In particular,
$[Q,Q],\,(Q,Q,Q),\,Z(Q),\,N(Q),\,C(Q)$ are normal subloops of 
an automorphic Moufang  loop $Q$.\\

\section{Main Theorem}\label{1/2-iso}

The aim of this note is to show the following

\begin{thm}\label{thm2}

Let $L$ be a  finite  automorphic Moufang loop  and  let $\tau$ be a
half-automorphism of $L$. Then 
$\tau$ is an automorphism of $L$ or $\tau$ is an anti-automorphism
of $L$.
\end{thm}

\begin{conj}\label{c1}
For every automorphic  Moufang loop all half-isomor\-phisms are trivial. 
\end{conj}
In order to prove this Theorem  we need to collect some results.

We will use the notation
$[u,v,w]=[[u,v]w]$ and by induction 
$$[a_1,...,a_{n+1}]=[[a_1,...,a_n],a_{n+1}]$$ 
for left-normed commutators.\\

\begin {lemma}\label{Br}{\rm (Bruck)} Let $L$ be a left  automorphic  Moufang loop and assume
$t,u,v \in L$.
Then  the following statements  hold:

{\rm (i)} If $L$ is generated by $3$ elements, then the associator subloop $(L,L,L)$  lies in the center of $L$.
\smallskip

{\rm (ii)}\,\,\,
$ [u,v]\in N(L)$.
\smallskip

{\rm (iii)}\,\,\,
$[uv,t]= [u,t][u,t,v][v,t]$,
\smallskip

{\rm (iv)}\,\,\,
$(nu,v,t)=(u,nv,t)=(u,v,nt)=(u,v,t)$  for all $n\in N(L)$.
\smallskip

{\rm (v)}\,\,\,
If $L$ is an automorphic Moufang loop, then
$ u^3\in N(L).$
\end {lemma}
\noindent
{\bf Proof.} {\cite{Bruck}, chapter  VII, Lemma 2.2
 \qed\\

Quite recently in the theory of finite Moufang loops a breakthrough has been reached:
the Theorem of Lagrange and an important part of the Theorems of Sylow have been proved 
(see \cite{GZLa}, \cite{GHLa}, \cite{Sylow}, \cite{p_Sylow}). In the following Lemma \ref{a} and Lemma \ref{dnil} these
facts are used. 

\begin {lemma}\label{a}
Let $L$ be a finite automorphic Moufang loop, let $\tau$ be a half-automorphism
of $L$ and let $S$ be a $3$-- Sylow subloop of $L$. Then
\smallskip

{\rm (i)}\,\,\,
$L = SN(L) =N(L)S$,

\smallskip
{\rm (ii)}\,\,\,
$\tau (L,L,L)=(\tau(L),\tau(L),\tau(L)),$
\smallskip

{\rm (iii)}\,\,
$\tau (N(L)) =N(\tau(L))$,
\smallskip

{\rm (iv)}\,\,\,
$\tau (C(L)) =C(\tau (L))$,
\smallskip

{\rm (v)}\,\,\,\,
$\tau [ L,L] = [\tau (L),\tau (L)].$
\end {lemma}
\noindent
{\bf Proof.} It is known by (v) from the
previous  Lemma \ref{Br} that any finite automorphic Moufang loop $L$
is the product of $S$ and the nucleus $N(L)$ of $L$.
Thus (i) holds.

Obviously  $\tau (S)$ is  a $3$--Sylow subloop of  $\tau (L).$ Applying (i) to 
$\tau(L)$ we obtain 
$$\tau (L)= \tau (S)\tau (N(L))
= \tau (S)N(\tau(L)).$$
For arbitrary elements $u,v,w\in L$ we put $u=s_1n_1,\,\,v=s_2n_2,\,\,w=s_3n_3$
with $s_i \in S,\,\,n_i \in N(L),\,\, 1\leq i\leq 3$.
Hence $(u,v,w) =(s_1,s_2,s_3)$ according to Lemma \ref{Br},(iv).
In  {\cite{GG} it was shown that for Moufang loop of odd order any half-automorphism is trivial,
so $\tau$ is trivial on $S.$\\
If $\tau$ is an isomorphism then
$$
\tau (u,v,w) =\tau(s_1,s_2,s_3)= (\tau(s_1),\tau(s_2),\tau(s_3))=
(\tau(u),\tau(v),\tau(w)).
$$
Otherwise straightforward  computation shows that
$$
\tau (u,v,w)=(\tau(u)^{-1}, \tau(v)^{-1},\tau(w)^{-1})^{-1}.
$$ 
Hence in both cases we have
$$
\tau (L,L,L) =(\tau (L),\tau (L),\tau (L))
\textnormal{ and }
\tau (N(L)) = N ( \tau (L))
$$
by the definition of $N(L).$

It is evident that
$\tau(Z(L))=Z(\tau (L))$
and $\tau (C(L)= C(\tau (L)).$\\
Finally, $\tau [L,L] =[\tau L, \tau L]$ since Moufang loops are diassociative.
\qed\\

In the Theory of Loops there exist different notions of nilpotency (\cite {Bruck},\cite {MP} and \cite {M}) which 
for groups, but not for loops in general, are equivalent.
Here we
will use the following  definition:\\
A Moufang loop $M$ is {\it commutatively nilpotent} if for some $n\in {\mathbb N}$ 
any left-normed  commutator word of length $n$ is trivial on $M$.

\begin{lemma}\label{dnil} 
Let $L$ be a finite commutatively nilpotent automorphic Moufang loop. Then\\
{\rm (i)} $L$  is the direct product of its $3$-Sylow subloop and its $3'$- Hall subgroup, which is contained in the nucleus of $L.$ \\
{\rm (ii)} $L$  is the direct product of its Sylow subloops.
\end{lemma}
\noindent
{\bf Proof.} First we show 
{\rm (i)}. By  Lemma \ref {a} for  every finite automorphic Moufang loop one has  $L = SN(L)$,
where $S$ is a 3-Sylow subloop of  $L$.
Let $a\in S$ and  $b\in L$ such that the orders of  $a$ and $b$ are coprime. The elements
$a$ and $b$ generate a  nilpotent subgroup of $L$.
Thus they commute. On the other hand, $N(L)$ is a nilpotent subgroup of $L$.
Hence $N(L)$ is the  direct product of its Sylow subgroups.
Let $\hat{N}$ be  the direct product of all Sylow subgroups of $N(L)$ apart of its $3$-Sylow subgroup. Then $L=S\hat{N}$ and  $S\cap \hat{N}=1$. Since $\hat{N}$ is subgroup of $N(L)$ it follows immediately that every element $a\in L$ has a unique presentation in the 
product $S\hat N$. 
Indeed  $\,s_1n_1=s_2n_2\,$ for $\,n_1, n_2\in \hat{N}$, and $s_1, s_2 \in S\,$  implies $\,s_2^{-1}s_1 =n_2n_1^{-1}.$ 

Using again the fact that $\hat{N}\leq N(L)$  by straightforward verification one can show 
that $S$ is invariant under the inner mappings group $\mathrm {Inn}(L)$ 
and $\hat{N}$  is  also invariant under the inner mappings group $\mathrm {Inn}(L).$ 
For example,
$$l_{a,b} (s) = (ab)^{-1}(a(b(s)))= (s_1n_1\cdot s_2n_2)^{-1}(s_1n_1(s_2n_2(s)))=s_3,$$
where $s,s_1, s_2,s_3\in S, \,n_1, n_2\in \hat{N}.$ 

It follows that
$L$  is the direct product of its $3$-Sylow subloop and its $3'$- Hall subgroup, which is contained in the nucleus of $L.$ 
By analogous arguments one can show that every Sylow subloop is normal in $L$, and thus we have the second statement of  the Lemma
(see \cite{Bruck},  p.72).
 \qed\\
By  Lemma \ref {a},(ii)   one can         
define the 
induced half-isomorphism 
$$
\bar{\tau}:L/(L,L,L) \to \tau(L)/(\tau(L),\tau(L),\tau(L))
$$
on a finite automorphic Moufang loop $L$ with 
a  half-isomorphism $\tau$.
Denote   $\bar{L}=  L/(L,L,L)$   and   $\bar{\tau(L)}=\tau(L)/(\tau(L),\tau(L),\tau(L))$

\begin{lemma}\label{centrcom} 
Let $L$ be a finite automorphic Moufang loop with 
a nontrivial half-isomorphism $\tau$. Consider the 
induced half-isomorphism 
$$
\bar{\tau}:\bar{L} \to \bar{\tau(L)}.
$$
\noindent
Then the following statements hold:\\
{\rm (i)}  $\bar{\tau}$ is a trivial half-isomorphism.\\[1ex]
 {\rm (ii)}  If $L$ is generated by $3$ elements and if $\bar{\tau}$ is an isomorphism on the group $\bar{L}$, then
for every pair $u,v\in L$ satisfying the condition 
$$\tau (uv)=\tau (v)\tau (u)$$
one has $[u,v]\in C(L)$ and $[\tau(u),\tau (v)]\in C(\tau(L))$.\\[1ex]
{\rm (iii)}  If $L$ is generated by $3$ elements and if $\bar{\tau} $ is an anti-isomorphism,
then $[w,t]\in C(L)$ and $[\tau(w),\tau (t)]\in C(\tau(L))$
for all $w,t\in L$
satisfying  the condition 
$$\tau(wt)=\tau(w)\tau(t)$$
\end{lemma}
\noindent 
{\bf Proof.} {\rm (i)}
Since $L/(L,L,L)$ is a group it follows from Scott's Theorem \cite{Scott} that
$\bar{\tau}$ is an isomorphism or an anti-isomorphism.\\
{\rm (ii)}Assume first that $\bar{\tau}$ is an isomorphism. For $u,v\in L$ and $\bar{u}=u(L,L,L),\\ 
\bar{v}=v(L,L,L)$ we get
$\bar{\tau} (\bar{u}\bar{v})=\bar{\tau}(\bar{u})\bar{\tau}(\bar{v})$. 
We have  on the other hand
$\bar{\tau} (\bar{u}\bar{v})=\bar{\tau}(\bar{v})\bar{\tau}(\bar{u})$.
Thus the commutator $[\bar{\tau}(\bar{u}),\bar{\tau}(\bar{v})]$  is trivial, and 
$[\tau(u),\tau(v)]$ is in the $(\tau(L), \tau(L),\tau(L))= \tau(L,L,L) $ and consequently
by Lemma \ref {Br},(i)  one knows that
$[\tau(u),\tau(v)]$ is in the center of $\tau (L)$.
Using half-isomorphism $(\tau )^{-1}$  analogously one can show that $[u,v]\in C(L)$. Hence 
statement {\rm (ii)} is shown.

As we mentioned above every  anti-isomorphism is the composition of some  isomorphism and inverse mapping, 
so one can treat  the case {\rm (iii)} of an anti-isomorphism $\bar{\tau}$ in a similar way.
\qed\\ 

In what follows we will consider the case that $\bar{\tau}$ is an isomorphism on $\bar{L}$.
If   $\bar{\tau}$ is an anti-isomorphism on $\bar{L}$ one can prove the Theorem \ref{thm2} using 
the inverse mapping and analogous arguments.

\begin{lemma}\label{d}
Let $L$ be a finite automorphic Moufang loop generated by 3 elements with 
a nontrivial half-isomorphism $\tau$.
Put $$D(L,\tau)=\{g\in L\mid\exists h\in L: \tau(gh)=\tau(h)\tau(g)\not=\tau(g)\tau(h)\}.$$
If $d\in [L,L],$ then $[d,D(L,\tau)]\subseteq C(L).$
\end{lemma}
\noindent
{\bf Proof.}
 Let $d\in [L,L],\, g\in D(L,\tau),$ choose $h\in L$ such that $\tau(gh)=\tau(h)\tau(g)\not=\tau(g)\tau(h).$
Since $[L,L]\subseteq N(L)$ a subloop $L_0$ generated by $d,h,g$ is a group. By Scott's Theorem
$\tau$ restricted to $L_0$ is an anti-monomorphism. Hence $\tau(dg)=\tau(g)\tau(d).$
If $\tau(g)\tau(d)\not=\tau(d)\tau(g),$ then by Lemma \ref{centrcom} we get $[g,d]\in C(L).$ Finally, if 
$\tau(g)\tau(d)=\tau(d)\tau(g),$ then $[g,d]=1.$
\qed
\\
\begin{lemma}\label{GG}{\rm(Gagola-Giuliani \cite{GG})}
Let $M$ be a finite  Moufang loop. Let $\tau$ be a  half-isomorphism
of $M.$ Then\\
{\rm (i)}
 $\tau$ is a semi-isomorphism of $M$, i.e. for any $u,v\in M$
 $\tau(uvu)= \tau (u)\tau (v)\tau (v)$\\
{\rm (ii)} 
If  $M$ is  of odd order, then  $\tau$ is trivial. In this case in particular
$\tau (M)$ is isomorphic to $M$.\\
{\rm (iii)} 
If  $\tau$ is not trivial, then there exist $x,y,z\in M$ such that $[x,y]\not=1,$ $[x,z]\not=1$ and

\begin{alignat}{1}
&\tau(xy)=\tau(x)\tau(y)\not=\tau(y)\tau(x),\\
&\tau(xz)=\tau(z)\tau(x)\not=\tau(x)\tau(z).
\end{alignat}
\end {lemma}
\qed
\\

Let us call such elements $\{x,y,z\}$  as it is described above in the Lemma \ref{GG} a
{\it Gagola-Giuliani triple} or $GG$-triple for short.\\
\begin{lemma}\label{ml}
Let $L$ be a finite automorphic Moufang loop with a non-trivial half-automorphism $\tau$. 
A subloop of $L$ generated by a  Gagola-Giuliani triple  is commutatively nilpotent.
\end{lemma}
\noindent
{\bf Proof.} 
Let  $L$ be a finite automorphic Moufang loop with non-trivial half-automorphism $\tau $ and let $\{x,y,z\}$ be a $GG$-triple of a loop $L$.
Let us recall that this means that conditions (iii)  of Lemma \ref{GG} hold.
Hence (1) implies that
$\tau|_{\left<x,y\right>}$ is a monomorphism and (2) implies that $\tau|_{\left<x,z\right>}$ is an anti-monomorphism. Thus there are two possibilities
for the subgroup $\left<y,z\right>$:  either  $\tau|_{\left<y,z\right>}$ is a monomorphism or
 $\tau|_{\left<y,z\right>}$ is an anti-monomorphism.

Put  $\tau(x)=a,\,\tau(y)=b,\,\tau(z)=c$.
Note that $\{a,b,c\}$ is a $GG$-triple if and only if  $\{x,y,z\}$ is a $GG$-triple. 
Indeed the conditions (1) and (2) of Lemma \ref{GG} are  equivalent to
$$\tau^{-1}(ab)=\tau^{-1}(a)\tau^{-1}(b)\not=\tau^{-1}(b)\tau^{-1}(a)$$
$$\tau^{-1}(ac)=\tau^{-1}(c)\tau^{-1}(a)\not=\tau^{-1}(a)\tau^{-1}(c)$$
Let $L_1$ be a loop generated $x,y,z$ and $L_2$ be a loop
generated by $a,b,c$.
Obviosly the mappings $\tau|_{L_1}:L_1\to L_2$  and
$\tau^{-1}|_{L_2}:L_2\to L_1$
are non-trivial half-isomorphisms.
\\[1ex]
\noindent
{\bf Case 1:} Suppose that $\tau|_{\left<y,z\right>}$ is an anti- monomorphism and
 $$\tau(yz)=\tau(z)\tau(y)=cb\not=\tau(y)\tau(z)=bc.$$
 One has
 $[L_2,L_2,u]\subseteq C(L_2),$
where 
 $u\in\{a,b,c \} \subseteq D(L_2,{\tau}^{-1})$ 
by Lemma \ref{d}. Hence  $[L_2,L_2,L_2]\subseteq C(L_2)$ by Lemma \ref{Br},(iii)
and $L_2$ is commutatively nilpotent of class $\leq 4.$ \\[2ex]
\noindent
{\bf Case 2:} Suppose that $\tau|_{\left<y,z\right>}$ is a monomorphism.
Let us study the behavior of the
half-isomorphism $\tau$ restricted to the following two-generated  subgroup:
 $$L_3=\left<xy,zx\right>,$$
\noindent
Suppose that  $\tau$ restricted to $L_3$ is a monomorphism. Then
$$\tau(xy\cdot zx)=\tau(xy)\tau(zx)=ab\cdot ac.$$
By the Moufang identity and by Lemma \ref{GG}(i), 
$$\tau(xy\cdot zx)=\tau(x\cdot yz\cdot x)=a\cdot bc\cdot a=ab\cdot ca.$$
Hence $ac=ca$ which contradicts  the condition that $\{x,y,z\}$ is GG- triple.
Therefore  $\tau$ restricted to $L_3$  is an anti-monomorphism. 
Now  by Lemma \ref {centrcom} one has 
$[ab,ac]\in C(L_2)$ as well as $[a,c]\in C(L_2).$ But it can also be proven that
$$[ab,ac]=[a,c]\not=1.$$
\noindent
In order to see this, note that
\begin{eqnarray*}
  ac\cdot ab &=& \tau (zx)\tau (xy)=\tau (xy\cdot zx)
 \\
 &=&\tau (x\cdot yz\cdot x)=a\cdot bc\cdot a=ab\cdot ca,
\end{eqnarray*} 
Let us consider the $GG$-triple
$\{a,ab,ac\}$. The map
$\tau ^{-1}$ restricted to $\left<a,ab\right>$ is a monomorphism,  $\tau ^{-1}$ restricted to $\left<a,ac\right>$ is an anti-monomorphism,
and $\tau ^{-1}$ restricted to $\left<ab,ac\right>$ is an anti-monomorphism too.
Thus we have the conditions of the  Case 1 and therefore  the loop generated by $\{a,ab,ac\}$ is commutatively nilpotent. 
But $$L_2= \left<a,ab,ac\right>= \left<a,b,c\right>.$$
The Lemma is proved.
\hfill\qed\\[1ex]
Now everything is ready to prove the Theorem \ref {thm2}.

By Lemma \ref {dnil} 
the loop $L_2$ is a direct  product of its Sylow subloops.
So  $\tau ^{-1}$ acts componentwise on  $L_2$. All Sylow subloops except the
 $3$-Sylow subloop are  groups, and the  $3$ -Sylow subloop is a loop of odd order.
By the Lemma \ref{GG}(ii) $\tau ^{-1}$ is a trivial half-isomorphism on $L_2$, which forms a contradiction with the conditions of Lemma \ref{ml}. 
Thus Theorem \ref{thm2}
 is proved.
\qed

\section{Examples.}
We give two examples of loops which admit a non-trivial half-auto\-mor\-phism. Both examples are in a different way not too far away from being automorphic Moufang loops. All our statements about these examples are easily checked using the LOOPS package of GAP \cite {L}\\[1ex]
\noindent
{\bf 1.} Let $Q_1$ be the loop of order 16 given by the Cayley table
\smallskip\noindent
\begin{center}
\begin{tabular}{|c ||c| c| c| c| c| c| c| c| c| c| c| c| c| c| c| c| c}
\hline
  $\circ$ & 1 & 2 & 3 & 4 & 5 & 6 & 7 & 8 & 9 & 10 & 11 & 12 & 13 & 14 & 15 & 16 \\
\hline\hline
    1 & 1 & 2 & 3 & 4 & 5 & 6 & 7 & 8 & 9 & 10 & 11 & 12 & 13 & 14 & 15 & 16 \\
\hline
        2 & 2 & 4 & 8 & 6 & 3 & 1 & 5 & 7 & 14 & 9 & 16 & 10 & 11 & 12 & 13 & 15 \\
\hline    
3 & 3 & 5 & 4 & 7 & 6 & 8 & 1 & 2 & 15 &13 & 9 & 11 & 14 & 16 & 12 & 10 \\
\hline    
4 & 4 &6& 7& 1& 8 &2 &3 &5 &12& 14 &15& 9& 16& 10 &11& 13\\
 \hline
    5 & 5&7& 2& 8& 4& 3& 6& 1& 13& 11& 14& 16& 12& 15& 10& 9\\
 \hline    
6 & 6 &1& 5& 2& 7& 4& 8& 3& 10& 12& 13& 14& 15& 9& 16& 11\\ 
\hline
7 & 7& 8& 1& 3& 2& 5& 4& 6& 11& 16& 12& 15& 10& 13& 9& 14\\
\hline    
8 & 8&3& 6& 5& 1& 7& 2& 4& 16& 15& 10& 13& 9& 11& 14& 12\\ 
\hline
9 & 9 & 10& 11& 12& 16& 14& 15& 13& 4& 6& 7& 1& 5& 2& 3& 8\\
\hline
10 & 10 & 12& 16& 14& 15& 9& 13& 11& 2& 4& 5& 6& 3& 1& 8& 7\\
\hline  
11 & 11 & 13& 12& 15& 10& 16& 9& 14& 3& 8& 4& 7& 6& 5& 1& 2\\
\hline    
12 & 12&14& 15& 9& 13& 10& 11& 16& 1& 2& 3& 4& 8& 6& 7& 5\\
\hline
13 & 13 &15& 10& 16& 9& 11& 14& 12& 8& 7& 2& 5& 4& 3& 6& 1\\ 
\hline
14 & 14 &9 &13 &10 &11 &12 &16 &15 &6 &1& 8& 2& 7& 4& 5& 3\\ 
\hline
15 & 15 & 16 &9 &11 &14 &13 &12 &10 &7 &5 &1& 3& 2& 8& 4& 6\\ 
\hline
16 & 16 & 11& 14& 13& 12& 15& 10& 9& 5& 3& 6& 8& 1 &7 &2 &4\\
\hline    
\end{tabular}
\end{center} 
\bigskip\noindent

Then $Q_1$ is a Code loop (see \cite {Griess}) in particular it is a left automorphic  Moufang loop.
One can see that the mapping $\phi: Q_1 \to Q_1$ defined by
$$\phi(5) = 8, \phi(8) = 5, \phi(x) = x\,\,{\rm for}\, x\notin\{5,8\}$$
is a half-automorphism. Since

$\phi(2\circ7) = \phi(5) = 8 = \phi(7)\circ\phi(2) \neq \phi(2)\circ \phi(7) = 5$

$\phi(3\circ9) = \phi(15) = 15 = \phi(3)\circ\phi(9) \neq  \phi(9)\circ\phi(3) = 11$\\
this half-automorphism is non-trivial.\\[1ex]
2. The loop $Q_2$ defined by the Cayley table \\[2ex]
\begin{center}
\begin{tabular}{|c|| c| c| c| c| c| c| c| c|}
\hline
$\circ$&1  &  2  &  3  &  4  &  5 &  6 &  7 &  8\\
\hline \hline
1&1  &  2  &  3  &  4  &  5 &  6 &  7 &  8\\
\hline
2&2  &  1  &  4  &  3  &  6 &  5 &  8 &  7\\
\hline
3&3  &  4  &  1  &  2  &  7 &  8 &  6 &  5\\ 
\hline
4&4  &  3  &  2  &  1  &  8 &  7 &  5 &  6\\
\hline
5&5  &  6  &  8  &  7  &  1 &  2 &  4 &  3\\
\hline
6&6  &  5  &  7  &  8  &  2 &  1 &  3 &  4\\
\hline
7&7  &  8  &  5  &  6  &  3 &  4 &  2 &  1\\
\hline
8&8  &  7  &  6  &  5  &  4 &  3 &  1 &  2\\
\hline
\end{tabular}
\end{center}
{\vskip 0.3cm}
is an automorphic loop, but not a Moufang loop and the permutation 
$\phi =(3,5)(4,6)(7,8)$ is a half-automorphism of $Q_2$. By

$\phi(4\circ 6) = \phi(7) = 8 = \phi(4)\circ\phi(6) \neq \phi(6)\circ \phi(4)=7$

$\phi(4\circ 8) = \phi(6) = 4 = \phi(8)\circ\phi(4) \neq \phi(4)\circ \phi(8)=3$
\\
$\phi$ is no-trivial half-automorphism.
{\vskip 0.8cm}
\noindent

{\bf Acknowledgements.}
The authors thank to Stephen Gagola III, Jacob Mostovoy and 
Peter Plaumann for discussions about  some details of the article
and to Alberto Elduque for his help to improve the article. 
The first author thanks to FAPESP, CNPQ and russian grant
RFFI 13-01-00239. The second author thanks 
 to SRE, Mexico to support  her visit to UAEM, Mexico during  May-July, 2012.
The third author thanks to the russian grant 2014/319. 
 The fourth author thanks for support to
2011-2013 UCMEXUS-CONACYT Collaborative Grant CN-11-567,\\
2012-2013 FAPESP Grant processo 2012/11068-2 and \\
2012-2013 CONACYT Grant for Sabbatical year at the Institute of
Mathematics and Statistics of the University of Sao Paulo, Brazil.


\begin{thebibliography}{999}
\bibitem{Bruck}{\sc R.~Bruck}, {\it A Survey of Binary Systems},  Springer-Verlag (1966)
\bibitem{CS}{\sc R.~Carrillo Catalan, L.~Sabinina}, On
smooth power alternative loops, Communications in Algebra {\bf 32},
no. 8 (2004), 2969-2976.
\bibitem{p_Sylow}{\sc S.~Gagola~III}, The development of Sylow $p$-subloops in finite Moufang loops, J. Algebra, {\bf 322}, no. 8 (2009), 2804--2810. 
\bibitem{GHLa}{\sc S.~ Gagola III, J.I.~Hall},
Lagrange's theorem for Moufang loops,  
Acta Sci. Math. (Szeged) {\bf 71} (2005), no. 1-2 (2005), 45-64.
\bibitem{GG}{\sc S.~ Gagola III, M.L.~Merlini Giuliani},
Half--isomorphism  of Moufang loops of odd order, J. of Alg. and Appl, {\bf 11}, no. 5 (2012), 194-199.
\bibitem{L}{\sc S.~ Gagola III, M.L.~Merlini Giuliani}, On half-automorphisms of certain Moufang loops with even order, J. of Algebra, {\bf 386} (2013),
131-141.
\bibitem{Griess}{\sc R.I.~Griess}, Code loops. J.of Algebra. {\bf 100}, no. 1 (1986), 224-234.
\bibitem{GPS}{\sc A.~Grishkov, P.~Plaumann L.~Sabinina},  Structure of free automorphic loops.
Proc. Amer. Math. Soc. {\bf 140}, no. 7 (2012), 2209- 2214.
\bibitem{GZLa}{\sc A.~Grishkov, A.~Zavarnitsine},
Lagrange's theorem for Moufang loops,
Math. Proc. Cambridge Philos. Soc. {\bf 139}, no. 1 (2005), 41-57.
\bibitem{Sylow}{\sc A.~Grishkov, A.~Zavarnitsine}, Sylow's theorem for Moufang loops.
J.of Algebra {\bf 321}, no. 7 (2009),  1813-1825.
\bibitem{Kikkawa}{\sc M.~Kikkawa}, Geometry of homogeneous Lie loops, Hiroshima Math J. {\bf 5} (1975), 141-179 
\bibitem{KKP}{\sc M.~Kinyon, K.~Kunen,J.~D.~Phillips},  Every  diassociative  $A$- loop is Mou\-fang. Proc. Amer. Math. Soc. {\bf 130} (2004), 619-624 .
\bibitem{M}{\sc J.~Mostovoy}, Nilpotency and Dimension Filtration for Loops, Comm.
Algebra {\bf 36} , no. 4 (2008), 1565–-1579.
\bibitem{MP}{\sc J.~Mostovoy, J.M.~Perez-Izquierdo}, Dimension Filtration on Loops, Israel
J. Math. {\bf 158} (2007), 105–-118.
\bibitem{L}{\sc G.~Nagy, P.~Vojtˇechovsk´y}, LOOPS.Version 2.2.0. Computing with quasigroups
and loops in GAP. preprint
\bibitem{PS}{\sc P.~Plaumann, L.~Sabinina}, On nuclearly nilpotent loops of finite exponent. Comm. in Algebra {\bf 36} (2008), 1346--1353.
\bibitem{Sabinin}{\sc L.~Sabinin}, {\it Smooth quasigroups and loops}. Kluwer (1999).

\bibitem{Scott}{\sc W.~R.~Scott}, Half-Homomorphisms of Groups, Proc. of the AMS {\bf 8}, no. 6 (1957), 1141--1144.
\end{thebibliography}
\end{document}